\documentclass[11pt,letterpaper,reqno]{amsart}

\makeatletter%
\@mparswitchfalse%
\makeatother%

\reversemarginpar%

\usepackage{amsfonts}
\usepackage{amssymb}
\usepackage{amsmath}
\usepackage{amsthm}
\usepackage{mathtools}
\usepackage{amstext}
\usepackage{tikz-cd}
\usepackage{enumitem}

\usepackage[left=1.3in,right=1.3in,top=1.2in, bottom=1.3in]{geometry}

\definecolor{refkey}{HTML}{0003d6}
\definecolor{labelkey}{HTML}{083b19}

\usepackage[disable,colorinlistoftodos, textwidth=28mm, color=green!30]{todonotes}

\usepackage[ruled,noend,linesnumbered,algosection]{algorithm2e}
\usepackage{graphicx}

\definecolor{cit}{HTML}{117733}
\definecolor{lin}{HTML}{0003d6}
\definecolor{gree}{HTML}{6dc066}
\usepackage{hyperref}
\hypersetup{bookmarksnumbered,colorlinks=true,linkcolor=lin,citecolor=cit,urlcolor=black}

\usepackage[section]{placeins}
\usepackage{caption}

\newcommand{\mycomment}[1]{} %
\newcommand{\mc}[1]{} %

\DeclareMathOperator{\Pic}{Pic}

\DeclareMathOperator{\PALP}{PALP}

\newcommand{\CC}{\mathbb{C}}
\newcommand{\Z}{\mathbb{Z}}

\newcommand{\R}{\mathbb{R}}
\newcommand{\PP}{\mathbb{P}}

\newcommand{\Q}{\mathbb{Q}}

\newcommand{\adj}{\leftrightarrow}
\newcommand{\nr}{\sim}
\DeclareMathOperator{\adjt}{adj}

\DeclareMathOperator{\nbd}{nbd}

\DeclareMathOperator{\AD}{AD}
\DeclareMathOperator{\conv}{conv}

\DeclareMathOperator{\Int}{int}

\newcommand{\fl}[1]{\lfloor#1\rfloor}

\newcommand{\vv}[1]{\lvert#1\rvert}
\newcommand{\vr}[1]{\langle#1\rangle}
\newcommand{\ra}{\to}

\newcommand{\txt}[1]{\text{ #1}} %
\newcommand{\sname}{Picard stage}

\theoremstyle{plain}%
\newtheorem{theorem}{Theorem}[section]

\newtheorem{corollary}[theorem]{Corollary}
\newtheorem{prop}[theorem]{Proposition}
\newtheorem{lemma}[theorem]{Lemma}
\newtheorem{conj}[theorem]{Conjecture}
\newtheorem{question}[theorem]{Question}

\theoremstyle{definition}%
\newtheorem{defi}[theorem]{Definition}

\newtheorem{remark}[theorem]{Remark}

\newtheorem{heuristic}[theorem]{Heuristic}

\newcommand{\pf}{{\em Proof}}%
\newcommand{\pfof}[1]{{\em Proof of #1}}%

\begin{document}
\title[Toric Fano pairs with relatively large sum of Picard ranks]{Toric Fano pairs with relatively large sum of Picard ranks}
\author{Zhuang He}
\email{he.zhu@northeastern.edu}

\begin{abstract}
We classify polar pairs of Gorenstein toric Fano varieties of dimension $d$ whose Picard rank sum equals the minimum of the number of vertices and facets of the associated reflexive polytopes minus $d$. In view of the recently discovered upper bound on the Picard rank sum, these are the pairs closest to a simplicial-simple pair. We show that except for dimension three, one of the pair must have Picard rank one, which is an analogue to Eikelberg's result on simple polytopes. We show more properties of the facet complexity of a convex polytope under various polytope products.
\end{abstract}
\maketitle

\section{Introduction}\label{Intro}
In this paper we study the Picard group and related invariants of toric Fano varieties. For complete toric varieties, their Picard rank equals to the second Betti number. The {\em facet complexity} $\rho'$ was introduced in \cite{He2026} as a numerical invariant of convex polytopes (details in Section \ref{fcPrelim}). For a reflexive polytope $P$, $\rho'_P$ is at least the Picard rank of the corresponding toric variety $X_P$ of $P$.

The attention on non-$\Q$-factorial Picard groups of toric Fano varieties dates back to Eikelberg \cite{Eikelberg1992, Eikelberg1993}. Batyrev conjectured that for smooth toric Fano varieties of dimension $d$, the Picard rank is at most $2d$. The result was first proved by Casagrande \cite{Casagrande2005} for all $\Q$-factorial toric Fano varieties. Progress on the Picard ranks includes \cite{Voskresenskiĭ1985, Debarre2000, Nill2005M, 2008Øbro, Fujita2019, He2026}, until the recent result \cite{HeNill2026} proved the $2d$ upper bound for all Gorenstein toric Fano varieties. The toric variety $X_P$ is Gorenstein and Fano if and only if $P$ is a reflexive polytope. As a comparison, the number of vertices of a reflexive polytope, which equals to the rank of the $\Q$-Weil divisor group of $X_P$, can grow quite fast (see \cite{HofscheierNill2026}).

This paper is a continuation of the aforementioned paper \cite{He2026}. We first explore more properties of the facet complexity with operations on polytopes. Theorems \ref{fcSum}, \ref{fcProduct}, and \ref{fcFreeJoin} exhibit identities for the facet complexity that are satisfied by the Picard ranks. With the additivity under free sums (Theorems \ref{fcSum}), we construct examples (Corollary \ref{prfcdiffer}) in all dimension $d\geq 3$ where the facet complexity differs from the Picard rank, answering a question in \cite{He2026}.

The main theorem of \cite{He2026} (Thm. 1.1) provides an upper-bound for the sum of Picard ranks of a polar pair $(X_P, Y_Q)$, where $(P, Q=P^*)$ is a polar pair of reflexive polytopes of dimension $d\geq 3$:
\begin{equation}\label{PQbound}
\rho_X + \rho_Y \le \min \{\vv{F(P)}, \vv{V(P)}\} -d +1,
\end{equation}
where the equality holds if and only if $P$ or $Q$ is simplicial. We can always assume that $\vv{F(P)}\le \vv{V(P)}$, so $P$ is simpler, and $Q$ is more simplicial. If we define the {\em \sname{}} of $(P,Q)$ as:
 \[\varepsilon(P,Q):= \vv{F(P)} -d +1 - \rho_X - \rho_Y.\]
 Then the theorem above says $\varepsilon(P,Q)$ is non-negative, and $\varepsilon(P,Q)=0$ if and only if the pair is simple-simplicial.
 
\begin{figure}[h]
    \centering
    \begin{minipage}{0.47\textwidth}
        \centering
        \includegraphics[width=0.43\textwidth]{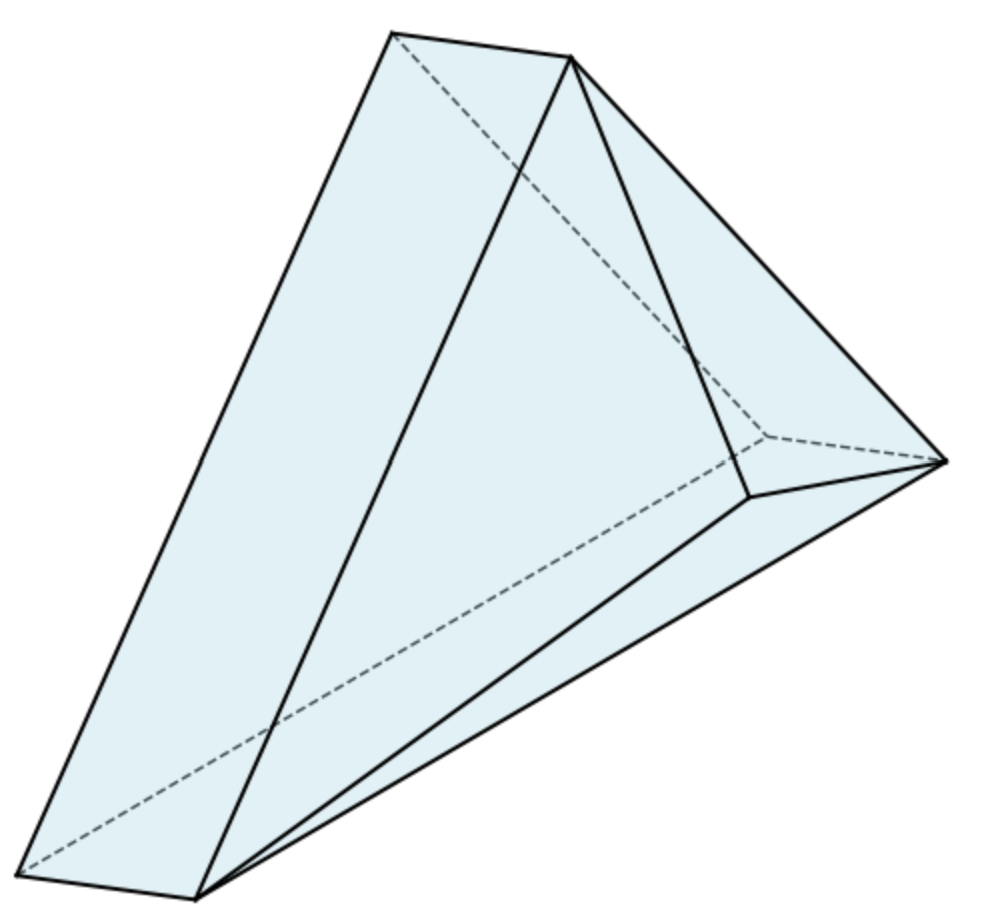}
    \end{minipage}\hfill
    \begin{minipage}{0.47\textwidth}
        \centering
        \includegraphics[width=0.43\textwidth]{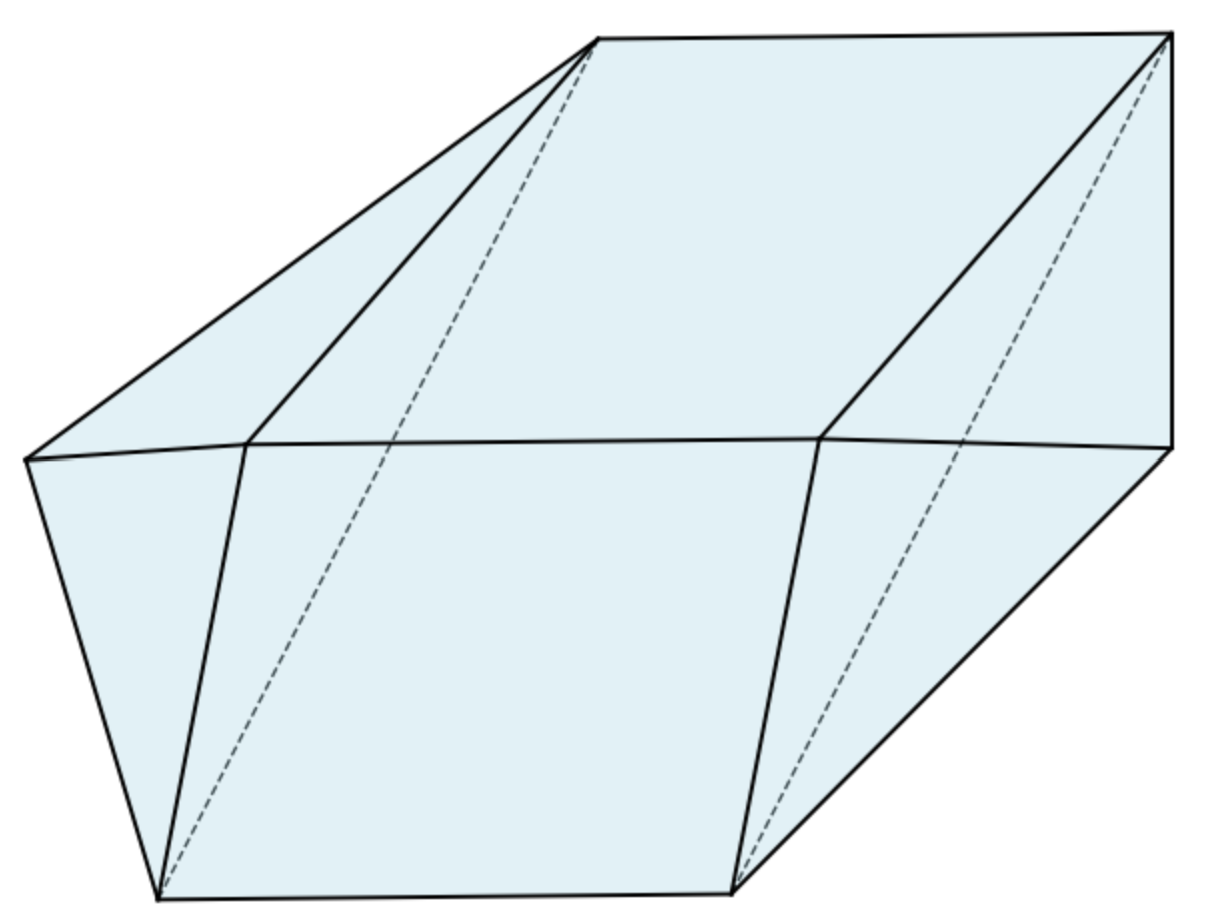}
    \end{minipage}
	\caption{Capped triangular prisms with $\varepsilon=1$. Left: PALP id 1145, with one cap; Right: PALP id 2840, with two caps.}
	\label{prismfigures}
\end{figure}

 It is a natural question to classify reflexive polytopes by their \sname{}. All reflexive polygons are simplicial, so $\varepsilon=0$. We characterize those pairs with $\varepsilon=1$:
 
\begin{theorem}\label{epsilon1}
 If $\varepsilon(P,Q)=1$, then every facet $F$ of $Q$ has exactly $d+1$ or $d+2$ vertices. If $d\ge 4$, then there are two cases: either $Q$ has exactly one non-simplex facet which is a free sum of two simplices (FS), or $Q$ has non-simplex facets that are all repeated pyramids (PS) over a common face $B$ which is a direct sum of two simplices. If $d=3$, then either $Q$ has exactly one non-triangle facet which is a quadrilateral, or the combinatorial type of $Q$ is a triangular prism capped with triangle pyramids on one or both ends (Figure \ref{prismfigures}).
\end{theorem}

We recall the theory in \cite{Eikelberg1992} of non-$\Q$-factorial toric Picard groups. The space $\AD(Y)$ is spanned by the affine relations that are linear on each maximal dimensional cone of the fan of $Y_Q$. Then $\rho_Y = \vv{V(Q)} - \dim \AD(Y) - d$. See \cite{Fujita2019, He2026} for details. With this, we immediately have a different interpretation of (\ref{PQbound}) that for any reflexive polytope $P$, and $X=X_P$:
\begin{equation}\label{PQboundAlt}
\rho_X  \le \dim \AD(P^*) +1,
\end{equation}
and $\varepsilon(P,Q)=\dim \AD(P^*)+1-\rho_X$. We note that $X_P$ is projective, so $\rho_X\ge 1$. Thus if $\dim \AD(P^*)=1$, then $\varepsilon=1$ too, as otherwise $P$ is simple. It is not hard to classify polytopes with $\dim \AD=1$ to see they contain precisely the non-prism cases in Theorem \ref{epsilon1} (See Section \ref{dplus2V} and Remark \ref{AD1}). Theorem \ref{epsilon1} can then be stated as:

\begin{theorem}\label{epsilon1withAD}
 If $\varepsilon(P,Q)=1$, then either $\dim \AD(Q)=1$, or $d=3$, $\dim \AD(Q)=2$, and the combinatorial type of $Q$ is a triangular prism capped with triangle pyramids on one or both ends.
\end{theorem}

Theorem \ref{epsilon1withAD} rules out the possibility of $\dim \AD(Q)\ge 3$ when $\varepsilon=1$. In particular, if $\varepsilon(P,Q)=1$, then $\rho_X=2$ for the triangular prisms, or $1$ otherwise. This phenomenon resembles the simple-simplicial case, where \cite{Eikelberg1993} first proved the general result that the toric variety of a simple polytope has Picard rank one. In general, one can ask how fast $\rho_X$ grows (within the upper bound $2d$) as $\varepsilon$ increases, while the classification are expected to be scattered for larger $\varepsilon$.

Structure of the paper: we review preliminaries on convex geometry and facet complexities in Section \ref{Prelim}. 
Section \ref{FCprops} studies properties of the facet complexity.
Sections \ref{efficient} to \ref{mainproof} classify Gorenstein toric Fano pairs with \sname{}  one.

{\bfseries Declaration on AI usage} 
Claude Opus 5 and 5.5 (high) were used to (1) generate the proof of Theorem \ref{fcSum}, (2) learn and verify the definitions, lemmas and references in Section \ref{prelimProducts}, and (3) prove the last two steps of the Proof of Theorem \ref{epsilon1} in dimension three (every non-simplex pair meet along an edge, and up to two vertices outside of the prism). The author learned how the redundancy of a facet can be applied in the context of facet complexities by reading ChatGPT-5.6 Sol's alternative proof of \cite[Lemma 4.8]{He2026}.

The complete paper was written solely by the human author. The main idea proving Theorem \ref{epsilon1} was found without AI assistance. All heuristics, questions and conjectures in the paper were formed by the author without AI hints.

{\bfseries Acknowledgment} The author thanks Cinzia Casagrande for discussions. The author thanks Benjamin Nill and the Otto-von-Guericke-Universität Magdeburg for hosting the research visit, during which the proof of Theorem \ref{epsilon1} was found over an A0-size paper whiteboard. The project heavily relies on the open-source software {\em Macaulay2} \cite{M2} with the \verb|polyhedra| package \cite{PolyhedraSource, PolyhedraArticle}; {\em SageMath} \cite{SageMath}; and the PALP database \cite{Kreuzer2004}. The author was partially funded by the European Union - Next Generation EU, Mission 4, Component 1, CUP D53D23005860006, within the project PRIN 2022L34E7W ``Moduli spaces and birational geometry''.

\section{Preliminaries}\label{Prelim}
\subsection{Convex geometry and Fano polytopes} 

We refer to \cite{Debarre2000,Nill2005M} for introductions on Fano polytopes. Here we list and define most notations.
\begin{description}[leftmargin=!,labelwidth=6em]
\item[Lattices] $N=\Z^d$ and $M=N^*$ is the dual lattice of $N$.
\item[Polytopes] $P,Q,R$ are convex polytopes in $N_\R$ or $M_\R$.
\item[Dimension] The dimension of a set $S$ in $N_\R$ is the dimension of the affine space spanned by $S$. A $d$-dimensional polytope is called a $d$-tope.
\item[Faces] We write $F,G,H,\cdots$ for facets, $u,v,w,x,y,\cdots$ for vertices. A ridge is named $\tau$.
\item[Sets] $F(P)$ is the set of facets, and $V(P)$ the set of vertices of $P$. 
\item[Adjacency] Two facets $F$ and $G$ are adjacent if they share the same ridge. Two vertices $v$ and $w$ are adjacent if they share an edge. We write $F\adj G$ and $v\adj w$. $\adjt(F)$ is the set of adjacent facets of $F$. 
\item[Neighbors] $\nbd(v)$ is the set of facets of $P$ containing a vertex $v$.
\item[Topology] $\conv(S)$ is the convex hull of $S$. $\partial S$ is the boundary of $S$, and $\Int S$ the interior of $S$.
\item[Polar] Assume $0\in\Int P$. The polar polytope of $P$ is $P^*:=\{y\in M_\R: \vr{x,y}\ge -1 \txt{for all } x \in P\}$.

\end{description}
\vspace{1em}

{\bfseries Toric varieties.} $X=X_P$ the projective toric variety whose fan is the spanning fan over $P$. 
For instance, let $P$ be the bipyramid with vertices $(e_1,e_2,-e_1-e_2,e_3,-e_3)$ in $N_\R\cong \R(e_1,e_2,e_3)$. Then $X_P\cong \PP^2\times \PP^1$.

{\bfseries Fano polytopes.} A lattice polytope $P$ with $0\in\Int P$ is Fano if all the vertices of $P$ are primitive lattice points.

{\bfseries Reflexive polytopes.} A lattice polytope $P$ with $0\in\Int P$ is reflexive if the polar polytope $P^*$ is a lattice polytope too. If $P$ is reflexive, then so is $P^*$; the only lattice point in $\Int P$ is $0$, so $P$ is particularly Fano. $X_P$ is Gorenstein and Fano if and only if $P$ is reflexive.

\subsection{Operators on polytopes} \label{prelimProducts}\,

{\bfseries Free sum.} Let $P\subset \R^{d}$ and $Q\subset \R^{e}$ be two full dimensional convex polytopes. Suppose $0\in \Int P$ and $0\in \Int Q$. The free sum (or direct sum) $P\oplus Q$ is the convex hull $\conv ( P\times \{0\}, \{0\} \times Q)\subset \R^{d+e}$. There are $\vv{F(P)}\cdot \vv{F(Q)}$ facets of $P\oplus Q$, consisting of $\{F \oplus G: F\in F(P)$, $G\in F(Q)\}$. The vertices of $P\oplus Q$ are exactly $(V(P)\times \{0\}) \cup ( \{0\} \times V(Q))$, hence $\vv{V(P\oplus Q)}=\vv{V(P)}+ \vv{V(Q)}$.

{\bfseries Cartesian product.} The Cartesian product (or direct product) $P\times Q$ is the set $\{(x,y)\in \R^{d}\times \R^{e}: x\in P, y\in Q\}$. Under the assumption that $0\in \Int P$ and $0\in \Int Q$, we have
\begin{align}\label{sumproddual}
(P\times Q)^*\cong P \oplus Q.
\end{align}
Thus $P\times Q$ has $\vv{F(P)}+ \vv{F(Q)}$ facets, all of the form $F\times Q$ or $P\times G$, for $ F\in F(P)$ and $G\in F(Q)$. As a direct result, if $P$ and $Q$ are both reflexive, then $P\oplus Q$ and $P\times Q$ are reflexive too.

{\bfseries Free join.} We define the free join $P*Q$ of level $(-1,1)$ by 
\[
P*Q := \conv ( P\times \{0\}\times \{-1\}, \{0\} \times Q \times \{1\})\subset \R^{d+e+1}.
\]
With this definition, $0\in \Int (P*Q)$ too. The facets of $P*Q$ are either $F*Q$ or $P*G$, where $F\in F(P)$, and $G\in F(Q)$. The vertices of $P*Q$ are $(V(P)\times \{0\}\times \{-1\}) \cup ( \{0\} \times V(Q)\times \{-1\})$. Thus the vertices of $P*Q$ are all primitive lattice points. In particular, $P*Q$ is a Fano polytope if both $P$ and $Q$ are lattice polytopes.

\begin{lemma}\label{dualFreeJoin}
Suppose $0\in \Int P$ and $0\in \Int Q$. Then $(P*Q)^*=(2P^*) * (2Q^*)$.
\end{lemma}

\pf.
A point $w=(u,v,t)\in (P*Q)^*$ if and only if for any $x\in P$ and $y\in Q^*$, $\vr{x,u}-t\ge -1$ and $\vr{y,v}+t \ge -1$. Let $x=y=0$. Then necessarily $-1\le t\le 1$. So the condition is equivalent to that $u\in (1-t) P^*$ and $v\in (1+t)Q^*$. Let $\lambda:=(1-t)/2$. There always exist $u',v'$ such that $u=\lambda u'$ and $v=(1-\lambda) v'$. Then
\[w=\lambda(u',0,-1)+(1-\lambda)(0,v',1). \]
so equivalently, $w$ is in the free join of $2P^*$ and $2Q^*$.\qed

By Lemma \ref{dualFreeJoin}, If  $P, Q, 2P^*$ and $2Q^*$ are all lattice polytopes, then $P*Q$ is a reflexive polytope.

\subsection{Facet complexities}\label{fcPrelim} Let $P$ be a full-dimensional convex polytope in $N_\R\cong \R^d$ with $m$ facets. A {\em facet sequence} over $U\subset F(P)$ is an ordering $L=(G_0,\cdots, G_{u-1})$ of $U$. Let $L=(G_0,\cdots, G_{m-1})$ be a facet sequence over $F(P)$. For any $i\geq 1$, let $w_i:=\dim (\cup_{j<i} V(G_j) )\cap V(G_i)$, with the convention that $\dim \emptyset= -1$. Let $a_i=d-1-w_i$ if $w_i\geq d-2$, and $\infty$ otherwise.  
If $G=L_i$, $i\geq 1$, then we also write $w(L,G):=w_i$, and $a(L, G):=a_i$. Let $S_i(L):=\cup_{j\le i} V(G_j)$, and $S(L):=S_{u-1}$. Then let $\rho'_L:=\sum_{i=1}^{m-1} a_i$. 

\begin{defi}
The {\em facet complexity} of $P$ is
\[\rho'(P):=\min \{\rho'_L: L \txt{a facet sequence over } F(P)\}.\]
\end{defi}

A facet sequence $L$ is called {\em adjacent} if for each $i\geq 1$, $G_i\adj G_j$ for some $j$, $0\le j<i$. An adjacent facet sequence over $F(P)$ always exists \cite[Lemma 3.13]{He2026}, which has a finite facet complexity because all the $a_i$ are either $0$ or $1$. Thus $\rho'_P$ is always a finite integer. We say $G$ is {\em redundant} with respect to $L$ (and $P$), if $a(L,G)=0$. Then $\rho'_L$ equals to $\vv{F(P)}-r(L)-1$, where $r(L)$ is the number of redundant facets in $L$. %

\begin{prop}\label{fcabc}
 If $P$ is a Fano polytope, then $\rho_X\le \rho'_P$.
\end{prop}
\pf. See \cite[Thm. 3.19]{He2026} for the proof in the reflexive case. If $P$ is a Fano polytope, then $X_P$ is $\Q$-Gorenstein, so Eikelberg's constructions also apply, by choosing the $\Q$-Cartier divisor $L$ as $-K$. Then the same proof shows the inequality.\qed

The facet complexity provides a good estimation and upper bound of the Picard ranks. By \cite[Example 3.23]{He2026}, among the $4,319$ three-dimensional reflexive polytopes. all except $22$ polytopes have the same Picard rank and facet complexity. The exceptional $22$ consist of $21$ with $\rho=1, \rho'=2$, and one with $\rho=2, \rho'=3$. Thus $\rho'_X\leq 6$ too for $3$-dimensional reflexive polytopes.

\subsection{Connectedness}
Let $P$ be any polytope. The edge graph of $P$ is connected. %
\begin{lemma}\cite[Lemma 3.13]{He2026}\label{outofNbd}
For any vertex $v$ of $P$, there is an adjacent sequence over $F(P)\backslash \nbd(v)$, 
\end{lemma}

\subsection{$d$-topes with $d+2$ vertices}\label{dplus2V}
It is a well-known result that there are exactly $\fl{d^2/4}$ combinatorial types of $d$-topes with $d+2$ vertices (see \cite[\S 6.1]{Gruenbaum2003}. They are the $k$-th repeated pyramids over a direct sum $B$ with combinatorial type $\Delta_m \oplus \Delta_n$, with $k\geq 0$, $m,n\geq 1$, and $k+m+n=d$. Then the $d+2$ vertices of such polytopes $P$ consist of $k$ apexes, $m+1$ vertices of $\Delta_m$, and $n+1$ vertices of $\Delta_n$. We say $B:=B(P)$ the base of $P$. In particular, If $k=0$, then $B(P)=P$.

We say $P$ is a DS if $k=0$, and PS-$k$ otherwise. $P$ is simplicial if and only if $k=0$.

\begin{lemma}\label{miss1inBase}
Let $d\ge 3$. Suppose $P$ is a $d$-tope with $d+2$ vertices. Let $v\in B(P)$. Then there are at least two distinct facets of $P$ that miss $v$.
\end{lemma}
\pf. If $P$ is a DS of type $\Delta_m \oplus \Delta_n$, then we can assume $v\in V(\Delta_m)$. Let $G$ be the unique facet of $\Delta_m$ not containing $v$, and $H$ be any facet of $\Delta_n$, then $G\oplus H$ is a facet of $P$ not containing $v$. There are $n+1\ge 2$ choices of $H$.
Suppose $P$ is a PS-$k$ over base $B$ with type $\Delta_m \oplus \Delta_n$, withe $k\ge 1$. Then there are $n+1$ facets $\tau$ of $B$ missing $v$ as above. Let $G$ be the $k$-th pyramid over $\tau$ with apexes those of $P$. Then $G$ is a facet of $P$ that misses $v$.\qed

\section{Facet complexities of products}\label{FCprops}
All polytopes named in $P,Q$ in this section are of positive dimensions.

\subsection{Free sums}
\begin{theorem}\label{fcSum}
Suppose $P$ and $Q$ are two convex polytopes such that $0\in \Int P$ and $0\in \Int Q$. Then $\rho'_{P\oplus Q} =\rho'_{P} +\rho'_{Q}$.
\end{theorem}

\pf. Suppose $d_1=\dim P$ and $d_2=\dim Q$. Let $L=(H_0,\dots,H_{v-1})$ be a facet sequence over $F(P\oplus Q)$. We define two induced facet sequences $L_P$ and $L_Q$ over $F(P)$ and $F(Q)$. Each $H_i=F_i\oplus G_i$ for some $F_i\in F(P)$ and $G_i\in F(Q)$. We take the sequence $L(P):=\{F_i\}_{i=0}^{v-1}$ and then let $L_P$ be the subsequence consisting of the first appearance of each unique facet. Then $L_P$ is a facet sequence over $F(P)$, since each facet of $P$ appear exactly once. Define $L(Q)$ and $L_Q$ in the same way.

We show that if $\rho'_L<\infty$, then $\rho'_L= \rho'_{L_P}+\rho'_{L_Q}$. We abuse notation and identify $V(P\oplus Q)$ with $V(P)\sqcup V(Q)$. For any $S\subseteq V(P\oplus Q)$, and any $F\in F(P)$ and $G\in F(Q)$, we can decompose $S$ into $S=S_P\sqcup S_Q$, and $V(F\oplus G))=V(F)\sqcup V(G)$. Then by Lemma \ref{affDimJoin}:
\begin{align*}
\dim (S\cap V(F\oplus G))&=\dim ((S_P\cap V(F)) \cup( S_Q\cap V(G)))\\
&=\dim (S_P\cap V(F))+ \dim ( S_Q\cap V(G))+1.
\end{align*}

Now suppose $H_i=F\oplus G$, and let $S=\cup_{j<i} V(H_j)$. Then $w(L, F\oplus G)=\dim (S\cap V(H_i))=\dim (S_P\cap V(F))+ \dim ( S_Q\cap V(G))+1$. Note that that $\dim(S_P\cap V(F))= w(L_p, F)$ if $H_i$ is the first facet in $L$ where $F$ appears, and $\dim(S_P\cap V(F))= d_1-1$ otherwise. By assumption, $\rho'_L<\infty$, so each $w(L, F\oplus G)\geq d_1+d_2-2$. Therefore $\dim(S_P\cap V(F))\ge d_1-2$. This means $a(L_P,F)<\infty$, so $a(L_Q,G)<\infty$ too. Finally,
\begin{align}\label{alfg}
a(L, F\oplus G)&=d_1+d_2-1-w(L, F\oplus G)\\
	&=d_1-1-\dim (S_P\cap V(F))+ d_2-1- \dim ( S_Q\cap V(G)).
\end{align}
Summing up $a(L, F\oplus G)$ and we get the claim.

Now we can prove the theorem.

(i) $``\ge"$: There exists a facet sequence $L$ over $F(P\oplus Q)$ such that $\rho'_{L}<\infty$. Then $\rho'_L= \rho'_{L_P}+\rho'_{L_Q}\ge \rho'_P+\rho'_Q$.

(ii) $``\le"$: Suppose $T=(F_0,\cdots, F_{m-1})$ and $V=(G_0,\cdots, G_{n-1})$ are facet sequences over $\vv{F(P)}$ and $\vv{F(Q)}$ such that $\rho'_P=\rho'_T$ and $\rho'_Q=\rho'_V$. Let $H_{i,j}:=F_i \oplus G_j$, and 
\begin{align*}
L:=&(H_{0,0}, \cdots, H_{0,j}, \cdots, H_{0,n-1},\cdots, H_{i,0}, \cdots, H_{i,j}, \cdots, H_{i,n-1},\cdots, \\
&H_{m-1,0}, \cdots, H_{m-1,j}, \cdots, H_{m-1,n-1}).
\end{align*}
Then $L$ is a facet sequence over $F(P\oplus Q)$.
By definition, we find $L_P=T$ and $L_Q=V$. Therefore $\rho'_L=\rho'_{V_P}+\rho'_{V_Q}=\rho'_T+\rho'_V=\rho'_P+\rho'_Q$. So we find $\rho'_{P\oplus Q}\le \rho'_P+\rho'_Q$.
\qed

\begin{lemma}\label{affDimJoin}%
Suppose $V=V_1\oplus V_2$, $A\subset V_1$ and $B\subset V_2$ are finite point sets. Then
\[
\dim A\oplus B = \begin{cases}
\dim A + \dim B, \quad\quad\txt{if } 0\in\conv A \txt{and } 0\in\conv B; \\
\dim A + \dim B + 1, \quad\txt{otherwise}.  \\
\end{cases}
\]
\end{lemma}
\pf. This is a direct application of \cite[Theorem 2.16 (1)]{Gallier2011}.\qed

\begin{remark}
Step (\ref{alfg}) in fact shows that if $\rho'_L$ is a finite number, then $a(L, F\oplus G)=1$ if and only if exactly one of $F$ and $G$ appears for the first time in the sequences $L(P)$ or $L(Q)$, and $a(L_P,F)=1$ or $a(L_Q,G)=1$ respectively.
\end{remark}

We can compare Theorem \ref{fcSum} to the Picard ranks of the product of toric varieties. Suppose $X=X_P$ and $Y=Y_Q$ are two projective toric varieties over $\CC$. Then the fan of $X\times Y$ is the product fan of $X$ and $Y$, and therefore the spanning fan over $P\oplus Q$. It is well-known that $\Pic(X\times Y)\cong \Pic(X)\times \Pic(Y)$, for example, by noticing $H^1(X,\mathcal{O}_X)=0$, and using \cite[Ex. III.12.6]{Hartshorne1977}. That is, the Picard ranks are additive under free sums.

\begin{figure}[h]
\includegraphics[width=0.25\textwidth]{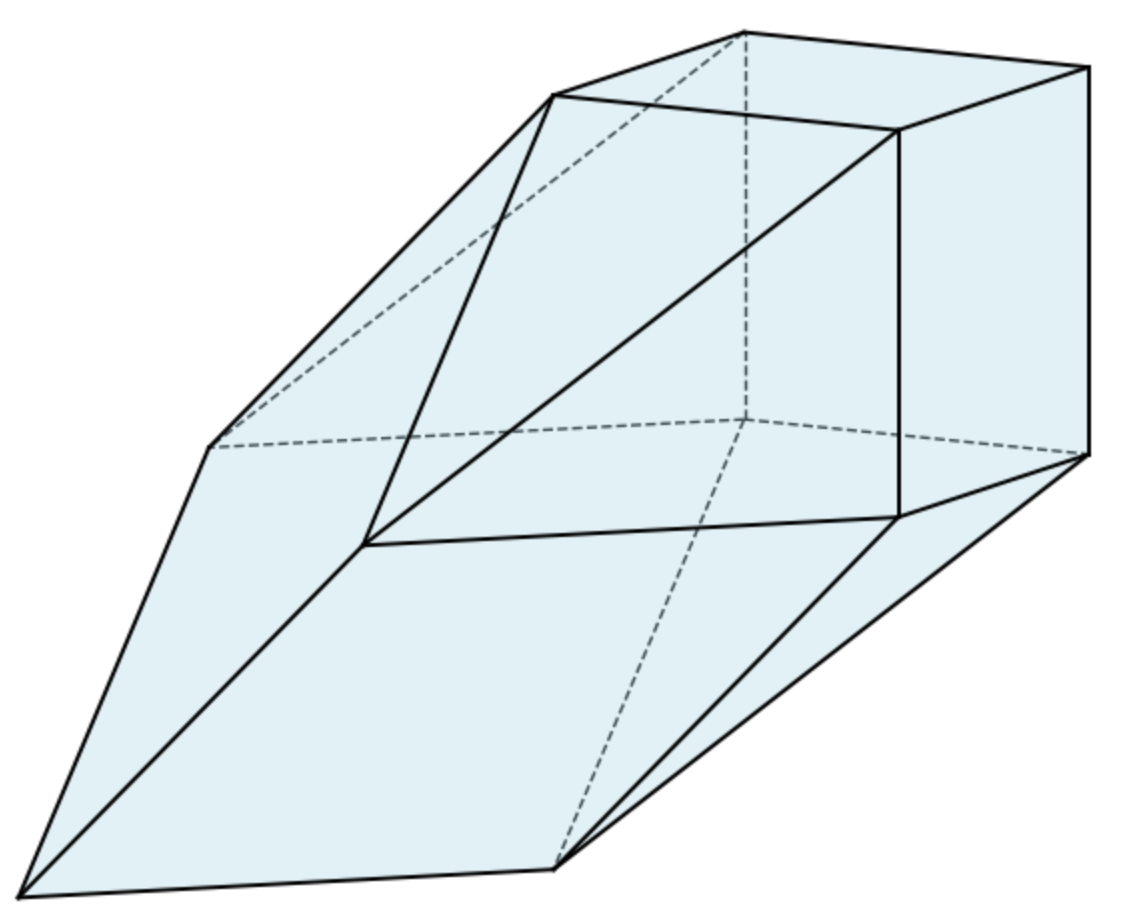}
\caption{The reflexive $3$-tope of PALP id 4252, with Picard rank $2$ and facet complexity $3$.}\label{4252figure}
\end{figure}

\begin{corollary}\label{prfcdiffer}
For any $d\geq 3$ there exists a reflexive $d$-tope $P$ such that $\rho_{X_P}\neq \rho'_P$.
\end{corollary}
\pf. By \cite[Example 3.23]{He2026}, there exist reflexive $3$-topes $P$ with $\rho_{X_P}<\rho'_P$. For example, take $P=\PALP(4252)$ as in Figure \ref{4252figure}. Then $\rho_{X_P}=2$ and $\rho'_P=3$. When $d>3$, let $Q:=P\oplus [-1,1]^{d-3}$. Then $Q$ is reflexive. By Theorem \ref{fcSum}, $\rho'_{Q}=\rho'_P+(d-3)\rho'_{[-1,1]}= 3+(d-3)\cdot 1=d$. The toric variety $X_Q\cong X_P\times (\PP^1)^{d-3}$. Therefore $\rho_{X_Q}=2+(d-3)\cdot 1=d-1<\rho'_Q$.\qed

\subsection{Cartesian products}
As a contrast, we show that the facet complexity of a Cartesian product is almost always one:

\begin{theorem}\label{fcProduct}
Suppose $P$ and $Q$ are two convex polytopes. Then $\rho'_{P\times Q} =1$ unless $\dim P=\dim Q=1$, in which case $\rho'_{P\times Q} =2$.
\end{theorem}
\pf. We have $F_1\times Q \adj F_2\times Q$ if and only if $F_1\adj F_2$, and $(F\times Q)\cap (P\times G)=F\times G$ is a ridge, so $F\times Q\adj P\times G$. Suppose $\dim Q\ge 1$. Then one can construct a facet sequence $L$ over $F(P\times Q)$ as follows: 
\[
L:=(F_0\times Q, F_1\times Q, P\times G_0,\cdots, P\times G_{n-1},F_2\times Q, \cdots, F_{m-1}\times Q),
\]
such that $F_0\adj F_1$. Then every facet in $L$ starting from $P\times G_0$ is adjacent to at least two previous facets of the opposite type, hence redundant. This gives $\rho'_L=1$, so $\rho'_{P\times Q}=1$.

The only case remaining is when both $\dim P=\dim Q=1$. Then $P\times Q$ is a quadrilateral, with $\rho'_{P\times Q}=2$.\qed

\begin{corollary}\label{prpolarPd}
Suppose $X_1,\cdots, X_n$ are Gorenstein toric Fano varieties, with $m\geq 2$. Let $Y$ be the polar variety of $\prod_{i=1}^n X_i$. Then $\rho_Y=1$, except for the case $m=2$ and $X_1\cong X_2\cong \PP^1$, in which case $\rho_Y=2$.
\end{corollary}

Corollary \ref{prpolarPd} can be seen as a generalization of Eikelberg's result \cite{Eikelberg1993} that if $P$ is a complete simple polytope, then $\rho_{X_P} = 1$. 

\begin{heuristic}
We say a polar pair of reflexive polytopes $(P,Q)$ is {\em quasi-simplicial-simple}, if $\rho'_P=1$ or $\rho'_Q=1$. Necessarily, the corresponding Picard rank equals one too. 
Non-uniform sampling of size about $1,000$ from Kreuzer and Skarke's database \cite{PALPdatabase, Kreuzer2000a} of reflexive $4$-topes gave only reflexive polar pairs with both Picard ranks one. This raises the question to further characterize such reflexive polar pairs.
\end{heuristic}
\begin{conj}\label{conjDensityAlmostSimple}
Let $A_d$ be the number of quasi-simplicial-simple polar pairs of reflexive polytopes, and $R_d$ be the total number of polar pairs of reflexive polytopes, both up to lattice isomorphism. Then 
\[
\lim_{d\ra\infty} \frac{A_d}{R_d}=1.
\]
\end{conj}

\subsection{Free joins}
\begin{theorem}\label{fcFreeJoin}
Suppose $P$ and $Q$ are two convex polytopes. Then $\rho'_{P* Q} =1$.
\end{theorem}
\pf. Choose any facet $F$ of $P$ and $G$ of $Q$. Let $L$ be a facet sequence over $F(P*Q)$ that starts with $(F*Q, P*G)$. Then after the second facet, all the vertices of $P*Q$ have been collected, so the remaining facets are all redundant. This gives $\rho'_L=1$, and hence $\rho'_{P* Q} =1$.\qed

With Lemma \ref{dualFreeJoin}, we have:
\begin{corollary}\label{fcFreeJoinDual}
Suppose $P$ and $Q$ are two Fano polytopes. Let $X=X_{P*Q}$, and $Y$ be the polar variety of $X$. Then $\rho_X=\rho_Y=1$.
\end{corollary}

\begin{heuristic}\label{bothRhoOne}
Theorem \ref{fcFreeJoin} generalizes \cite[Corollary 4.10]{Eikelberg1993} that the toric variety of the free join of complete polytopes has Picard rank one, in the case of Fano polytopes. We do not know a characterization of reflexive polytopes that can be written as a free join.
\end{heuristic}

\subsection{Discussion}
\begin{heuristic}
The three properties above of the facet complexity are surprising to the author, as the facet complexity was defined as an approximation to the Picard rank, with the goal to simplify calculations. In spite of the non-equality between $\rho$ and $\rho'$, the shared properties under products prompt the question whether the facet complexity can be associated with some topological invariants of toric Fano varieties.
\end{heuristic}
\begin{question}
Are there other nontrivial numerical invariants $f$ on convex polytopes such that for all convex polytopes $P, Q$ with the origin in their interior, the following hold?
\begin{enumerate}
\item $f(P\oplus Q)=f(P)+f(Q)$;
\item $f(P\times Q)=1$ unless $\dim P=\dim Q=1$; and 
\item $f(P* Q)=1$.
\end{enumerate}
\end{question}
\section{$d$-polytopes with at most $d$ redundant vertices}\label{efficient}
Consider a vertex sequence $T=(v_0,\cdots, v_s)$ over $V(P)$. We can define the redundancy of $v_i$ in $T$ by the dual sequence $T^*=(F_0,\cdots, F_s)$. That is, $v_i$ ($i\ge 1$) is redundant if $a(T^*, F_i)=0$.
\begin{defi}\label{vR}
Suppose $P$ is a $d$-tope, with $d\geq 2$. Let the {\em (vertex-)redundancy} $r(P)$ of $P$ be the maximal number of redundant vertices of any adjacent vertex sequence over $V(P)$. 
A $d$-tope $P$ is called {\em efficient}, if $r(P)\le d$.
\end{defi}
Clearly,
\begin{align}\label{rfcCum}
r(P)+\rho'_{P^*}\le \vv{V(P)}-1.
\end{align}
By \cite[Lemma 4.8]{He2026}, for any $d$-tope $P$, $d\geq 2$, there exists an adjacent sequence $T$ over $V(P)$ with at least $(d-1)$ redundant vertices. So $d-1\le r(P)\le \vv{V(P)}-2$. This implies that any $d$-tope with at most $d+2$ vertices is efficient. It is easy to see any polygon is efficient, where any adjacent vertex sequence has exactly one redundant vertex which is the last vertex.

We will see next a $d$-tope with redundancy $d-1$ is a $d$-simplex. The natural next step is to understand efficient $d$-topes.

\begin{lemma}\label{efficientSimplicial}
Let $d\geq 3$. A simplicial $d$-tope of $n$ vertices has redundancy $n-2$. In particular, a simplicial $d$-tope is efficient if and only if it is a simplex or a direct sum $\Delta_a \oplus \Delta_b$ of an $a$-simplex and a $b$-simplex.
\end{lemma}
\pf. Any vertex sequence on a $d$-simplex is adjacent, and every vertex starting from the third one is on an edge with all the previous vertices, therefore redundant. This shows a $d$-simplex has redundancy $d-1$. Now suppose $P$ is a simplicial $d$-tope. There exists an adjacent facet sequence $L=(H_0,\cdots, H_{m-1})$ over $F(P)$. Since $H_0$ is a simplex, there is an adjacent vertex sequence $T_0$ over $V(H_0)$ such that every vertex starting from the third is redundant. This $T_0$ is then also adjacent when considered over $P$. Next, each $H_i$ adds to $S_i$ with at most one new vertex $v_i$, which is adjacent to all the other vertices in $H_i$. Since $d\geq 3$, this makes $v_i$ redundant. One then expand $T_0$ to $T$ over $V(P)$, by sequentially adding those $v_i$. Thus $T$ is an adjacent sequence. As a result, $r(P)\geq n-2$. Since $r(P)\le n-2$, we have the equality.

If $P$ is efficient, then $r(P)\le d$, so $n\le d+2$. So $P$ is $d$-simplex, a direct sum $\Delta_a \oplus \Delta_b$ with $a+b=d$, or a $k$-th repeated pyramid over a direct sum $\Delta_a \oplus \Delta_b$, with $a+b=d-k$, and $k\geq 1$. The last case is not simplicial. \qed

\begin{lemma}\label{efficientH}
Let $d\geq 3$. If a $d$-tope $F$ is efficient, and $H$ is a facet of $F$, then either $H$ is a $(d-1)$-simplex, or $H$ has exactly $d+1$ vertices, and every vertex $v$ of $H$ forms an edge with exactly one vertex $w$ of $P$ not in $H$.
\end{lemma}
\pf. Suppose $H$ has at least $d+2$ vertices. The edge graph of $ V(F)\backslash V(H)$ is connected (\cite{Brøndsted1983}, see \cite[Lem 2.2]{He2026}). Thus there exists an adjacent sequence $T$ over $V(F)$ by arranging vertices not in $H$ first, and followed by a sequence $T_H$ over $V(H)$ such that $T_H$ is an adjacent sequence when considered on $H$. Clearly, any vertex $v$ of $V(H)$ except the first is on an edge with a vertex not in $H$, and a previous vertex in $H$. Therefore $v$ is redundant. Since this gives at least $d+1$ redundant vertices, we reached a contradiction.

Now $H$ has $d$ or $d+1$ vertices. If $H$ has $d$ vertices, then $H$ is a $(d-1)$-simplex. If $H$ has $d+1$ vertices, with a vertex $v$ on two edges with two different vertices outside of $H$. Then we can form the adjacent sequence $T$ in the same way above, while requiring that $v$ is the first vertex of $T_H$. This makes $v$ redundant too, and $T$ has at least $d+1$ redundant vertices, contradiction.\qed

\begin{prop}\label{efficientF}
A $d$-tope $P$ ($d\geq 3$) is efficient if and only if $\vv{V(F)}\le d+2$ when $d\ge 4$, or $P$ is the triangular prism for $d=3$.
\end{prop}
\pf. First consider the case $d\geq 4$. Suppose $P$ is efficient. By Lemma \ref{efficientH}, any facet $H$ of $P$ is either a $(d-1)$-simplex, or of exactly $d+1$ vertices. Suppose it is the latter case. Then $H$ is either simplicial, or a repeated pyramid over a simplicial polytope. Since $\dim H\ge 3$, there is an adjacent vertex sequence $T_H$ over $V(H)$ with $d-1$ redundant vertices. 
We extend $T_H$ to an adjacent vertex sequence $T$ over $V(P)$. Suppose there are at least two different vertices $v$ and $w$ not in $H$. Then without loss of generality, we can assume $v$ and $w$ are the last two vertices of $T$. Since $d\geq 3$, both $v$ and $w$ are on an edge with two other vertices. Therefore $v$ and $w$ are both redundant, which implies that $P$ has at least $d+1$ redundant vertices, a contradiction. Therefore, $P$ is a pyramid over $H$, and $\vv{P}=d+2$.

The remaining possibility is that $P$ is simplicial. In this case, by Lemma \ref{efficientSimplicial}, $P$ is a simplex or a free sum of two simplices.

When $d=3$, Lemma \ref{efficientH} says that each facet $H$ of $P$ has at most $4$ vertices. If $P$ has a quadrilateral facet $H$, then $H$ allows one redundant vertex. We claim that any other facet $G$ of $P$ must meet $H$. Suppose on the contrary that $G$ is a facet disjoint from $H$. Then we can form an adjacent vertex sequence $T$ over $F(P)$ such that $T$ starts with $V(G)$ and ends with $V(H)$. Then $T$ has at least $1+3$ redundant vertices, one from $G$ and three from the last three of $H$, by the same argument as in the first part of the proof to Lemma \ref{efficientH}. This is a contradiction.

By Lemma \ref{efficientH}, every vertex $v$ of $H$ is of degree $3$, so no facet meets $H$ at exactly a vertex. Together, we find there are exactly four facets of $P$ except $H$. They each meets $H$ along an edge. Therefore there are at most two vertices outside of $H$, so $P$ is either the triangular prism or the quadrilateral pyramid. The remaining case is when $F$ is simplicial, and by Lemma \ref{efficientSimplicial}, $F$ is the $3$-simplex or the triangle bipyramid.

Finally, if $\vv{V(P)}\le d+2$, then $P$ is efficient. The triangular prism allows three redundant vertices, and its dual is a triangle bipyramid, with $\rho'=2$. Therefore by (\ref{rfcCum}), the triangular prism is efficient too.\qed

\section{Facet complexity for pairs of \sname{} one}\label{mainproof}

\subsection{Dual sequences}
In this section, $P$ is a reflexive $d$-tope with $d\geq 3$. 
We generalize the construction used the proof of \cite[Theorem 1.1]{He2026}. For any vertex sequence $T=(v_0,\cdots, v_{t-1})$ over $V(P)$, we say a facet sequence $L$ over $F(P)$ is a {\em companion sequence} of $T$, if the following all hold:
\begin{enumerate}
\item $L$ is the concatenation of sequences $(L_0, L_1,\cdots, L_{t-1})$;
\item for each $i=0,\cdots, t-1$, $L_i$ consists of the facets $M_i:=\nbd(v_i)\backslash \cup_{j<i} \nbd(v_j)$.
\item $L_0$ has at least $(d-2)$-redundant facets.
\item $\rho'_{L, P}< \infty$.
\end{enumerate}

\begin{lemma}\label{companionExists}
Any vertex sequence $T$ over $V(P)$ such that $\rho'_{T^*, Q}< \infty$ has an adjacent companion sequence.
\end{lemma}

\pf. Let $A_{n}$ be the concatenation of $(L_0, L_1,\cdots, L_{n})$. By \cite[Prop. 4.9]{He2026}, for any vertex $v_0$ of $P$, there exists an adjacent sequence $L_0$ over $\nbd(v_0)$ with at least $(d-2)$-redundant facets. Inductively, suppose we have arranged the sequences $L_0,\cdots, L_{n-1}$ such that $A_{n-1}$ is adjacent. For $n$, since $\rho'_{T^*, Q}\le \infty$, $v_n\sim^G v_{k_n}$ for a $k_n< n$, and a facet $G$ of $P$. Therefore $G\in M_s$ for some $s\le k_n$, and $G\in \nbd(v_n)$. With this, one can arrange $\nbd(v_n)$ into an adjacent sequence starting with $G$, and then remove those not in $M_i$. Let $L_n$ be the sequence of remaining facets in $\nbd(v_n)$. Then each $H\in L_n$ is adjacent to a previous facet which is either in some $L_i$, $i<n$, or in $L_n$. Thus $A_n$ is adjacent too. By induction, we have $L=A_{t-1}$ is adjacent. Hence $\rho'_{L, P}< \infty$, so (1)-(4) all hold.\qed

\begin{lemma}\label{companionSum}
For any vertex sequence $T$ over $V(P)$ such that $\rho'_{T^*, Q}\le \infty$, and any adjacent companion sequence $L$ of $T$, $\rho'_{T^*, Q} + \rho'_{L, P} \le \vv{F(P)}-d+1$.
\end{lemma}

\pf. The proof is almost the same as the proof in \cite[\S 4.3]{He2026}. The only difference is that we assume $\rho'_{T^*, Q}\le \infty$ instead of that $T$ is adjacent, which is stronger. The only step where this could make a difference is when $\vv{L_i}=0$ for $i>0$. We prove $v_i$ is redundant when $\vv{L_i}=0$, which will finish the proof. To see this, pass to $Q$ and consider $T^*=(F_0,\cdots, F_{t-1})$. Since $\vv{L_i}=0$, all the vertices of $F_i$ are added to some $L_j^*$, where $j<i$. By the definition of facet complexity, $F_i$ is redundant since its contribution $a_i=0$.\qed

\subsection{Pairs with \sname{} one}
In the following, $(P,Q)$ is a pair of reflexive polytopes of dimension $d\geq 3$. $X=X_P$ and $Y=Y_Q$ are their associated Gorenstein toric Fano varieties. We assume $\vv{V(Q)}\le \vv{F(Q)}$, so $Q$ is the ``more simplicial" one. Thus $\varepsilon(P,Q) = \vv{V(Q)} -d +1 - \rho_X - \rho_Y$. We define $\varepsilon'(P,Q) = \vv{V(Q)} -d +1 - \rho'_X - \rho'_Y$.
\begin{lemma}\label{1implies1}
If $\varepsilon=1$, then $\varepsilon'=1$, $\rho_X=\rho'_X$, and $\rho_Y=\rho'_Y$.
\end{lemma}

\pf. If $\varepsilon=1$, then they proof of \cite[Thm. 1.1]{He2026} shows that $P$ is not simple, and there exists a facet sequence $T$ over $F(Q)$ and a facet sequence $L$ over $F(P)$ such that $\rho_{T^*}+\rho_L< \vv{F(P)}-d+1$. In particular, since $ \vv{F(P)} -d=\rho_X+\rho_Y\le \rho'_X+\rho'_Y\le \rho_{T^*}+\rho_L$, it must be that $\rho'_X+\rho'_Y =\vv{F(P)} -d$, and $\varepsilon'=1$ too. Thus $\rho_X=\rho'_X$, and $\rho_Y=\rho'_Y$.\qed

\subsection{Proof of the classification in \sname{} one: the case $d\ge 4$}
Now we can prove the main classification.

\pfof{Theorem  \ref{epsilon1} when $d\ge 4$}. Suppose $\varepsilon(P,Q)=1$. Then $Q$ is not simplicial. If $Q$ has an non-efficient facet $F$, then $F$ has a vertex sequence of $d$ redundant points. Choose any adjacent facet sequence $T^*=(F,\cdots)$ and any companion sequence $L$. Then $L_0$ contains at least $d$ redundant vertices, hence $\varepsilon\geq d-(d-2)=2$, a contradiction. As a result, all facets of $Q$ are efficient, so by Proposition \ref{efficientF}, they are either a $(d-1)$-simplex, a direct sum FS, a PS-$k$ with $k\ge 1$, or a triangular prism when $d=4$.

{\bfseries Case I: no triangular prism facets.} We assume first $Q$ has no triangular prism facets, and show that either $Q$ has a unique non-simplex facet which is a FS, or all the non-simplex facets of $Q$ are of type PS-$k$ sharing a common base.

Suppose towards contradiction that this is not the case. Then by Section \ref{dplus2V}, there are two non-simplex facets $G$ and $H$, such that $B(G)\ne B(H)$.  There there exists a vertex $v\in B(H)\backslash B(G)$.

By Lemma \ref{outofNbd}, there exists an adjacent sequence $T_1^*$ over $F(Q)\backslash \nbd(v)$. By the choice of $G,H$ and $v$, $G$ is in $T_1^*$ but $H$ is not. We can then require additionally that $G$ be the first facet of $T_1^*$. We then extend $T_1^*$ into an adjacent sequence over $F(Q)$, by first attaching $H$, and then everything else. This $T$ is adjacent at $H$. To see this, note that because $d\geq 4$, and $v\in B(H)$, there are at least two facets $\tau_1,\tau_2$ of $H$ that miss $v$. Then there are two facets $F_i$ of $Q$ meeting $H$ along the ridge $\tau_i$, $i=1,2$. Since $v\not\in F_i$, $F_i$ is in $T_1^*$, so $T^*$ is adjacent at $H$, and moreover, $H$ is redundant in $T^*$.

Now $\rho'_{T^*, Q}<\infty$. Applying Lemma \ref{companionExists} and Lemma \ref{companionSum}, there exists a companion sequence $L$ of $T^*$, and $\rho'_{T^*, Q} + \rho'_{L, P} \le \vv{V(Q)}-d+1$. We can arrange the vertices in $L_0$ such that $L_0$ has $d-1$ redundant vertices by Proposition \ref{efficientF}, which is one less than the upper bound $d-2$. Suppose $H=T_{k}$. On the one hand, $H$ is redundant. On the other hand, $L_{k}$ contains $v$, and hence nonempty. This means the contribution of $L_{k}$ to  $\rho'_{T^*, Q} + \rho'_{L, P}$ is also one less than the upper bound. Together, we find $\varepsilon(P,Q)\geq 2$, which is a contradiction.

{\bfseries Case II: $d=4$ with a triangular prism.} In this case, suppose $Q$ has a triangular prism facet $H$ and another non-simplex facet $G$. Then $G$ is a triangular prism, triangle bipyramid or quadrilateral pyramid. We can construct then with the same strategy above an adjacent sequence $T*$ which starts from $G$, covers all the facets of $Q$ not contains a vertex $v\in H$, and is followed by $H$ and the rest. Since any $v\in V(H)$ misses exactly two facets of $H$, the same argument above shows that the companion sequence $L$ of $T^*$ contribute two less to $\rho'_{T^*, Q} + \rho'_{L, P}$, so we obtain $\varepsilon(P,Q)\geq 2$, a contradiction.

{\bfseries The other direction} Suppose $Q$ is of the declared shape. Then by \cite[Thm. 1.1]{He2026}, $\varepsilon\geq 1$. We show that $\varepsilon=1$. 

We compute the affine dependency by \cite[Definition 3.2]{He2026}. If $Q$ has a unique FS facet $H$ that is a free sum $\Delta_m \oplus\Delta_n$. There is a unique affine relations among the $d+1$ vertices of $H$. Any other facet $G$ of $Q$ is a simplex, so $\AD(G)=0$. Thus $\AD(X)=\AD(H)$ is of dimension one. By (\ref{PQboundAlt}), $\rho_Y\le 1$. Since $Y$ is projective, $\rho_Y=1$. So $\varepsilon=1$.

Otherwise, $Q$ has $s$ PS-$k$ facets $H_0,\cdots, H_{s-1}$, sharing the same base $B$. Then is a unique affine relation $s$ among $V(B)$ as above, so $\dim \AD(X)=1$, and $\varepsilon=1$ too.\qed

\subsection{Proof in $d=3$.}

\begin{lemma}\label{Jordan}
Suppose $A$ is a union of polygon regions along their edges, and $A$ is homeomorphic to the disk $D^2$. Let $a,b,c,d$ be four distinct points on $\partial A$ in sequence. Suppose there are two polygons $G,H$ in $A$ such that $\{a,c\}\subset G$ and $\{b,d\}\subset H$. Then $G=H$.
\end{lemma}
\pf.  Suppose $G\neq H$. Then $\overline{ac}$ and $\overline{bd}$ are not edges of $G$ and $H$, so their interiors must meet at some point $x$ in $\Int A$, by the Jordan curve theorem. Since $x\in\Int G$ and $x\in\Int H$, we have a contradiction.\qed

\pfof{Theorem  \ref{epsilon1} when $d= 3$}.  
Assume $\varepsilon(P,Q)=1$.

We first show that if there are two non-triangle facets $G$, $H$ of $Q$, then every pair of vertices $u\in V(G)$ and $v\in V(H)$ must be in a common facet of $Q$. Indeed, suppose $u\not\nr v$, then there is an adjacent sequence $T_1^*$ over $\nbd(u)$, with $G$ the last term. Then by Lemma \ref{outofNbd}, $L'$ can completed into into a facet sequence $T_2^*$ over $F(Q)\backslash \nbd(v)$, as no facets in $T_1^*$ contains $v$. Lastly, $S(T_2^*)$ contains exactly the other three vertices of $H$. Thus we can attach $H$ to $T_2^*$ and complete it into an adjacent sequence $T^*$ over $F(Q)$. Then both $G$ and $H$ are redundant in $T^*$ and both give one vertex in the companion sequence $L$ of $T^*$, so $\varepsilon\geq 2$, contradiction.

We now prove that every pair of non-triangle facets $G$, $H$ must meet along an edge. With the property shown above, suppose $G\cap H=\{v\}$ or empty, then we can cut $\partial Q$ along an edge path linking $G$ and $H$, removing $G$ and $H$, to get a disk $A$, whose boundary contains the vertices $(x_1,\cdots, x_m)$ of $G$ and $(y_1,\cdots, y_n)$ of $H$ in clockwise order, with $m,n\ge 4$, where we can assume $v=x_1=y_1$ if $v$ exists. Now choose any quadruple $(x_i,x_j,y_k,y_l)$ with $2\le i <j$, $2\le k <l$. Then $x_i\sim y_k$ and $x_j\sim y_l$ by the previous step. Thus by Lemma \ref{Jordan}, $x_i,x_j,y_k$ and $y_l$ are on a same facet $F$. Run this for every such quadruple, we find $F$ contains all the vertices $x_2,\cdots, x_m, y_2,\cdots, y_n$, which is a contradiction when $m\ge 3$ and $n\ge 3$, as $F\cap G$ and $F\cap H$ have at most two vertices.

Now we can assume $x_1=y_2$, $x_2=y_1$, so $G\cap H=\overline{x_1 x_2}$ is an edge. Remove $G\cup H$ from $\partial Q$ and apply Lemma \ref{Jordan} to the remaining vertices of $G$ and $H$. Then there is a facet $F$ containing all the other vertices of $G$ and $H$. This can only happen when $m=n=2$. Thus $G$ and $H$ are quadrilaterals. Furthermore, the same argument for $(G,F)$ shows that $F$ is a quadrilateral too. so $G,H,F$ form an uncapped triangular prism. Clearly, there is no other quadrilateral facets besides $G,H,F$, otherwise they must form another triangular prism, which is impossible.

Finally we show there is at most one additional vertex at each cap of the prism. The three facet $F,G$ and $H$ separate the rest vertices into two connected components. Choose one vertex out of each cap, say $v_1$ and $v_2$. Then $v_1\not\sim v_2$. Suppose $v_1\neq -v_2$, then by  by \cite[Proposition 3.1]{Nill2005M}, there is a lattice point $z\in \partial Q$ such that $z\sim v_1$ and $z\sim v_2$, which contradicts to that $\adjt(v_1)\cap \adjt(v_2)=\emptyset$. Thus $v_1=-v_2$. Since this holds for any $v_1, v_2$ out of each cap, they are both unique in their cap, so $\vv{V(Q)}\leq 8$. If neither $v_1, v_2$ exists, so $Q$ is a triangular prism, then $Q$ is simple, and $\varepsilon=0$.\qed

\begin{remark}\label{AD1}
If $d\ge 4$, then it is easy to see that $\dim \AD(Q)=1$ if and only if all the non-simplex facets of $Q$ are of $d+1$ vertices and share the same base, as obtained in Theorem \ref{epsilon1}. Indeed in the proof we have shown the sufficiently. Conversely, if any facet $F$ of $Q$ has $m\ge d+2$ vertices, then $\dim \AD(F)=m-d> 1$. If $Q$ has two facets that are DS, PS that do not share the same base, then they gave linearly independent affine relations, and $\dim \AD(Q)\ge 2$.
\end{remark}
\vspace{2em}

\bibliography{mybib.bib}

\begin{thebibliography}{{The}25}

\bibitem[BK]{PolyhedraSource}
Ren{\'e} Birkner and Lars Kastner.
\newblock {Polyhedra: convex polyhedra. Version~1.10}.
\newblock A \emph{Macaulay2} package available at
  \url{https://github.com/Macaulay2/M2/tree/master/M2/Macaulay2/packages}.

\bibitem[BK09]{PolyhedraArticle}
Ren{\'e} Birkner and Lars Kastner.
\newblock {Polyhedra: a package for computations with convex polyhedral
  objects}.
\newblock {\em The Journal of Software for Algebra and Geometry: Macaulay2}, 1,
  2009.

\bibitem[Br{\o}83]{Brøndsted1983}
Arne Br{\o}ndsted.
\newblock {\em Convex Polytopes}.
\newblock Springer New York, New York, NY, 1983.

\bibitem[Cas06]{Casagrande2005}
Cinzia Casagrande.
\newblock The number of vertices of a {{F}ano} polytope.
\newblock {\em Annales de l'Institut Fourier}, 56(1):121--130, 2006.

\bibitem[Deb01]{Debarre2000}
Olivier Debarre.
\newblock {\em Toric {F}ano varieties}, pages 93--132.
\newblock Bolyai Society Mathematical Studies 12. Springer, 2001.

\bibitem[Eik92]{Eikelberg1992}
Markus Eikelberg.
\newblock The {P}icard group of a compact toric variety.
\newblock {\em Results in Mathematics}, 22:509--527, 1992.

\bibitem[Eik93]{Eikelberg1993}
Markus Eikelberg.
\newblock {P}icard groups of compact toric varieties and combinatorial classes
  of fans.
\newblock {\em Results in Mathematics}, 23:251--293, 1993.

\bibitem[Fuj19]{Fujita2019}
Kento Fujita.
\newblock The generalized {M}ukai conjecture for toric log {F}ano pairs.
\newblock {\em European Journal of Mathematics}, 5:858--871, 2019.

\bibitem[Gal11]{Gallier2011}
Jean Gallier.
\newblock {\em Geometric Methods and Applications: For Computer Science and
  Engineering}, volume~38 of {\em Texts in Applied Mathematics}.
\newblock Springer New York, second edition, 2011.

\bibitem[Gr{\"u}03]{Gruenbaum2003}
Branko Gr{\"u}nbaum.
\newblock {\em Convex Polytopes}, volume 221 of {\em Graduate Texts in
  Mathematics}.
\newblock Springer, New York, second edition, 2003.

\bibitem[GS]{M2}
Daniel~R. Grayson and Michael~E. Stillman.
\newblock Macaulay2, a software system for research in algebraic geometry.
\newblock Available at \url{http://www2.macaulay2.com}.

\bibitem[Har77]{Hartshorne1977}
Robin Hartshorne.
\newblock {\em Algebraic Geometry}, volume~52 of {\em Graduate Texts in
  Mathematics}.
\newblock Springer-Verlag, New York, NY, 1977.

\bibitem[He26]{He2026}
Zhuang He.
\newblock Reflexive polytopes and the {P}icard ranks of {G}orenstein toric
  {F}ano varieties.
\newblock Accepted to {\em Épijournal de Géométrie Algébrique}. To appear,
  2026.

\bibitem[HN26a]{HeNill2026}
Zhuang He and Benjamin Nill.
\newblock {G}orenstein toric {F}ano $d$-varieties have {P}icard ranks at most
  $2d$, 2026.
\newblock arXiv preprint. arXiv:2611.21074.

\bibitem[HN26b]{HofscheierNill2026}
Johannes Hofscheier and Benjamin Nill.
\newblock Reflexive polytopes can have doubly exponentially many vertices,
  2026.
\newblock arXiv preprint. arXiv:2611.08819.

\bibitem[KS]{PALPdatabase}
Maximilian Kreuzer and Harald Skarke.
\newblock {PALP: a Package for Analyzing Lattice Polytopes}.
\newblock available at
  \url{http://hep.itp.tuwien.ac.at/~kreuzer/CY/CYpalp.html}.

\bibitem[KS00]{Kreuzer2000a}
Maximilian Kreuzer and Harald Skarke.
\newblock Complete classification of reflexive polyhedra in four dimensions.
\newblock {\em Adv. Theor. Math. Phys.}, (4):1209--1231, 2000.

\bibitem[KS04]{Kreuzer2004}
Maximilian Kreuzer and Harald Skarke.
\newblock {PALP}: A package for analysing lattice polytopes with applications
  to toric geometry.
\newblock {\em Computer Physics Communications}, 157(1):87--106, 2004.

\bibitem[Nil05]{Nill2005M}
Benjamin Nill.
\newblock Gorenstein toric {F}ano varieties.
\newblock {\em Manuscripta {M}athematica}, 116(2):183--210, Feb 2005.

\bibitem[{\O}br08]{2008Øbro}
Mikkel {\O}bro.
\newblock Classification of terminal simplicial reflexive {$d$}-polytopes with
  {$3d-1$} vertices.
\newblock {\em Manuscripta {M}athematica}, 125(1):69--79, Jan 2008.

\bibitem[{The}25]{SageMath}
{The Sage Developers}.
\newblock {\em {S}ageMath, the {S}age {M}athematics {S}oftware {S}ystem
  ({V}ersion 10.6)}, 2025.
\newblock {\tt https://www.sagemath.org}.

\bibitem[VK85]{Voskresenskiĭ1985}
Valentin~E. Voskresenskiĭ and Alexander~A. Klyachko.
\newblock Toroidal {F}ano varieties and root systems.
\newblock {\em Mathematics of the USSR-Izvestiya}, 24(2):221, Apr 1985.

\end{thebibliography}
\bibliographystyle{alpha}
\end{document}